\documentclass[graybox]{svmult}

\usepackage{type1cm}        
\usepackage{makeidx}         
\usepackage{graphicx}        
\usepackage{multicol}        
\usepackage[bottom]{footmisc}
\usepackage{color}
\usepackage{caption}
\usepackage{placeins}
\usepackage{newtxtext}       %
\usepackage{newtxmath}       

\makeindex             

\begin{document}

\title*{Inexact subdomain solves using deflated GMRES for Helmholtz problems}
\author{N. Bootland, V. Dwarka, P. Jolivet, V. Dolean and C. Vuik}
\institute{Niall Bootland \at University of Strathclyde, Department of Mathematics and Statistics, Glasgow, UK\\ \email{niall.bootland@strath.ac.uk}
\and Vandana Dwarka \at Delft University of Technology, Delft Institute of Applied Mathematics, Delft, The Netherlands\\ \email{v.n.s.r.dwarka@tudelft.nl}
\and Pierre Jolivet \at University of Toulouse, CNRS, IRIT, Toulouse, France\\ \email{pierre.jolivet@enseeiht.fr}
\and Victorita Dolean \at University of Strathclyde, Department of Mathematics and Statistics, Glasgow, UK\\ Universit\'e C\^ote d'Azur, CNRS, Laboratoire J.A. Dieudonn\'e, Nice, France\\ \email{work@victoritadolean.com}
\and Cornelis Vuik \at Delft University of Technology, Delft Institute of Applied Mathematics, Delft, The Netherlands\\ \email{c.vuik@tudelft.nl}}
%
%
\maketitle

\abstract{We examine the use of a two-level deflation preconditioner combined with GMRES to locally solve the subdomain systems arising from applying domain decomposition methods to Helmholtz problems. Our results show that the direct solution method can be replaced with an iterative approach. This will be particularly important when solving large 3D high-frequency problems as subdomain problems can be too large for direct inversion or otherwise become inefficient. We additionally show that, even with a relatively low tolerance, inexact solution of the subdomain systems does not lead to a drastic increase in the number of outer iterations. As a result, it is promising that a combination of a two-level domain decomposition preconditioner with inexact subdomain solves could provide more economical and memory efficient numerical solutions to large-scale Helmholtz problems.}

\section{Introduction}
\label{bootland_mini_03_sec:1}

In recent years, domain decomposition based preconditioners have become popular tools to solve the Helmholtz equation. Notorious for causing a variety of convergence issues, the Helmholtz equation remains a challenging PDE to solve numerically. Even for simple model problems, the resulting linear system after discretisation becomes indefinite and tailored iterative solvers are required to obtain the numerical solution efficiently. At the same time, the mesh must be kept fine enough in order to prevent numerical dispersion `polluting' the solution \cite{deraemaeker1999dispersion}. This leads to very large linear systems, further amplifying the need to develop economical solver methodologies.

Domain decomposition (DD) techniques combined with Krylov solvers provide a way to deal with these large systems \cite{MR3908314}. While the use of two-level deflation and DD techniques have been explored before, their expedience has primarily been measured in terms of providing a way to add a coarse space to obtain a two-level DD preconditioner \cite{MR4266496,MR4299039}. Without a coarse space, DD methods typically do not scale with the number of subdomains. Moreover, subdomain sizes need to be relatively small in order to optimally use local memory and direct solution methods on subdomains.

In this work we focus on the subdomain solves. Instead of using a direct solution method, we solve the local systems using GMRES preconditioned by a two-level deflation approach. As a result, similar to the inclusion of a coarse space on the fine-level, we obtain a two-level method on each subdomain as well. The inexact solve on the subdomains will allow for larger subdomains by reducing computing and memory requirements, especially in 3D. In order to allow for inexact subdomain solves we require a flexible wrapper for the outer iteration---the application of the DD preconditioner---and so we will use FGMRES. Our local subdomain solves will then employ a preconditioned GMRES iteration. It is well known that the cost of GMRES increases with each iteration. Thus, in order to mitigate the number of iterations at the subdomain level, we use a two-level deflation preconditioner \cite{MR3534874,MR4080796}.

The techniques proposed here will feature as the topic of future research on large-scale 3D applications using pollution-free meshes. In this work we introduce the key ideas and begin an initial exploration by considering a simple 2D model problem.

\section{Model problem, discretisation and preconditioning strategies}
\label{bootland_mini_03_sec:2}

Our model problem consists of the Helmholtz equation posed on the unit square:
\begin{subequations}
\begin{align}
-\Delta u - k^2u &= f & & \text{in } \Omega = (0,1)^{2}, \\
u &= 0 & & \text{on } \partial\Omega.
\end{align}
\label{bootland_mini_03_eq:01}
\end{subequations}
Here, the parameter $k$ denotes the wave number. Problem \eqref{bootland_mini_03_eq:01} is well-posed so long as $k^2$ is not a Dirichlet eigenvalue of the corresponding Laplace problem. Solving the problem with Dirichlet conditions provides a more robust test for the solver, as there is no shift keeping the spectrum away from the origin \cite{bootland2021analysis,MR3534874}.
In this work we will assume the problem, and any sub-problems, are well-posed. To discretise \eqref{bootland_mini_03_eq:01} we use piecewise linear (P1) finite elements on a uniform grid with mesh spacing given by $h = 2\pi k^{-1} n_\text{ppwl}^{-1}$, where $n_\text{ppwl}$ is the number of (grid) points per wavelength (hereinafter referred to as ``ppwl''). To test the solver performance we initially ensure $n_\text{ppwl} \approx 10$. We then double this to approximately 20 ppwl to obtain more accurate numerical solutions.
Letting $V^{h} \subset H_0^1(\Omega)$ denote the space of piecewise linear functions on our finite element mesh $T^h$ of  $\Omega$, our discrete solution $u_{h} \in V^{h}$ satisfies the weak formulation $a(u_h,v_h) = F(v_h) \ \forall \, v \in V^{h}$, where
\begin{align}
a(u,v) &= \int\limits_{\Omega} \left( \nabla u \cdot \nabla {v} - k^2 u {v} \right) \mathrm{d}\mathbf{x} & \text{and} & & F(v) &= \int\limits_{\Omega} f {v} \, \mathrm{d}\mathbf{x}. \label{bootland_mini_03_eq:02}
\end{align}
With the standard basis for $V^{h}$, we can write the weak formulation as finding the solution to the linear system $A\mathbf{u} = \mathbf{f}$. We now consider how to solve such systems.

The global matrix $A$ is preconditioned by a one-level domain decomposition method. To construct the decomposition, we define an overlapping partition ${\left\lbrace \Omega_{j} \right\rbrace}_{j=1}^{N}$ of $\Omega$ together with a restriction operator $R_j$ to move from the global level to the subdomain level. We only consider Cartesian (rectangular) subdomains. Using this decomposition, the restricted additive Schwarz (RAS) preconditioner is defined by
\begin{align}
M_{RAS}^{-1} = \sum_{j = 1}^{N} R_{j}^{T}D_{j}A_{j}^{-1}R_{j}, \label{bootland_mini_03_eq:03}
\end{align}
where $D_j$ are diagonal matrices representing a partition of unity ($\sum_{j = 1}^{N} R_{j}^{T} D_j R_j = I$) and $A_j = R_j A R_{j}^{T}$ are the local Dirichlet matrices. Each subdomain solve requires the solution of a local auxiliary linear system, which we denote by $A_j \tilde{\mathbf{u}}_j = \tilde{\mathbf{f}}_j$ as a general case. We solve these systems using a two-level deflation approach, that is, deflation is used to accelerate the convergence of GMRES by removing the near-zero eigenvalues. For normal matrices it has been shown that convergence can be directly related to the behaviour of these near-zero eigenvalues \cite{MR3296704}.

The two-level deflation preconditioner is defined as a projection operator $P_j$ which leads to solving $P_j A_j \tilde{\mathbf{u}}_j = P_j\tilde{\mathbf{f}}_j$, where
\begin{align}
P_{j} &= I - A_j Q_j & \text{with} & & Q_j &= Z_j E_{j}^{-1} Z_{j}^{T} & \text{and} & & E_j &= Z_{j}^{T} A_j Z_j. \label{bootland_mini_03_eq:04}
\end{align}
The rectangular matrix $Z_j$ in this particular setting is called the deflation matrix and its columns span the deflation space. The choice of $Z_j$ strongly dictates the overall convergence behaviour. Here, we use quadratic rational B\'ezier curves, as they have been shown to provide satisfactory convergence \cite{MR4080796}. Consequently, if we let $\tilde{u}^{i}_{j}$ represent the $i$-th degree of freedom (DOF) on subdomain $\Omega_{j}$, then in 1D $Z_j$ maps these nodal approximations onto their coarse-grid counterpart as follows
\begin{align}
\left[ Z_{j} \tilde{u}_{j} \right]^{i} = \frac{1}{8} \left( \tilde{u}_{j}^{2i - 2} + 4 \, \tilde{u}_{j}^{2i - 1} + 6 \, \tilde{u}_{j}^{2i} + 4 \, \tilde{u}_{j}^{2i + 1} + \tilde{u}_{j}^{2i + 2} \right).
\end{align}
As such, $Z_j$ can be constructed using the following 1D stencil $\frac{1}{8}\begin{bmatrix}
1 & 4 & 6 & 4 & 1
\end{bmatrix}$. The dimension of $Z_j$ will then be $n_j \times \frac{n_j}{2}$, where $n_j$ is the size of the local 1D system.

In 2D on Cartesian grids this can be naturally extended by using the Kronecker product. For a rectangular subdomain, letting $Z_j^{x}$ denote the 1D deflation matrix in the $x$-direction and $Z_j^{y}$ that in the $y$-direction, the 2D deflation matrix is given by
\begin{align}
Z_j = Z_j^{y} \otimes Z_j^{x}, \label{bootland_mini_03_eq:05}
\end{align}
assuming lexicographic ordering running through $x$ coordinates first.

\section{Numerical results}
\label{bootland_mini_03_sec:3}

We now provide numerical results for our model problem on the unit square. We take the right-hand side $f$ to be given by a point-source at the centre of the domain. Unless stated otherwise, 10 ppwl are used to construct the mesh and we let $n_{\text{glob}}$ be the number of DOFs along each edge of the square. For the domain decomposition, we use a uniform decomposition into $N$ (square) subdomains. Overlap is added by appending one layer of mesh elements in a Cartesian manner (note that this means subdomains touching only one edge of $\Omega$ are rectangular rather than square).

For the outer solve we use preconditioned FGMRES with the one-level RAS preconditioner \eqref{bootland_mini_03_eq:03}. The tolerance for the relative residual has been set at $10^{-6}$. For the inner solve on the subdomain level we use preconditioned GMRES with the two-level deflation preconditioner \eqref{bootland_mini_03_eq:04} instead of a direct solver. Note that subdomains systems are decoupled and so can be solved locally in parallel. We will vary the inner tolerance for the relative residual between $10^{-10}$ and $10^{-2}$ in order to assess an appropriate level of accuracy needed when solving the subdomain problems. The solver is equipped to deal with both symmetric (Dirichlet) and non-symmetric systems (Sommerfeld), as we use the RAS preconditioner together with GMRES, which do not require symmetry.

All matrices are constructed using FreeFem (http://freefem.org/) while the solvers are then implemented using PETSc (http://www.mcs.anl.gov/petsc/). Computations are carried out on a laptop with an i7-10850H processor having 6 cores (12 threads).

\subsection{Direct subdomain solves}
\label{bootland_mini_03_sec:3a}
We start by constructing a benchmark where we use a direct solution method for the subdomain solves, namely via an $LU$-decomposition. Table \ref{bootland_mini_03_tab:01} shows that the number of iterations does not scale as the number of subdomains $N$ increases, in agreement with the literature. The number of iterations also rapidly increases with the wave number $k$. The inclusion of a coarse space to improve both $k$-- and $N$--scalability on top of inexact subdomain solves will be explored in future research. An interesting observation is that increasing $n_\text{ppwl}$ leads to a higher iteration count. The opposite effect has been observed when using a two-level deflation preconditioner \cite{MR3534874,MR4080796}. There, a finer mesh leads to a smaller number of iterations as the mapping of the eigenvectors from the fine- and coarse-grid becomes more accurate.
\begin{table}[t]
\centering
\caption{FGMRES iteration counts using the one-level RAS preconditioner with direct subdomain solves.}
\begin{minipage}{.49\linewidth}
\centering
10 ppwl\\
\scalebox{1}{
\begin{tabular}{@{}cc|cccc@{}}
\hline
& & \multicolumn{4}{c}{$N$} \\ \hline \hline
$k$ & $n_{\text{glob}}$ & 4 & 9 & 16 & 25 \\ \hline
20 & 30 & 20 & 27 & 48 & 45 \\
40 & 60 & 31 & 60 & 85 & 101 \\
80 & 120 & 64 & 133 & 191 & 216 \\
160 & 240 & 159 & 262 & 365 & 495 \\ \hline
\end{tabular}}\label{bootland_mini_03_tab:01a}
\end{minipage}%
\begin{minipage}{.49\linewidth}
\centering
20 ppwl\\
\scalebox{1}{
\begin{tabular}{@{}cc|cccc@{}}
\hline
& & \multicolumn{4}{c}{$N$} \\ \hline \hline
$k$ & $n_{\text{glob}}$ & 4 & 9 & 16 & 25 \\ \hline
20 & 60 & 20 & 40 & 42 & 59 \\
40 & 120 & 37 & 66 & 89 & 115 \\
80 & 240 & 76 & 131 & 189 & 255 \\
160 & 480 & 130 & 289 & 398 & 520 \\ \hline
\end{tabular}}\label{bootland_mini_03_tab:01b}
\end{minipage}\label{bootland_mini_03_tab:01}
\end{table}

\subsection{Inexact subdomain solves}
\label{bootland_mini_03_sec:3b}
The direct subdomain solves remain feasible for medium-size problems. Once we move to high-frequency 3D problems the subdomain systems become larger and the direct solver will start to become inefficient and consume more computing power and memory. In order to assess the feasibility of inexact solves, we will use the benchmarks from Section \ref{bootland_mini_03_sec:3a} and compare with our iterative method. The aim of these experiments is twofold. First, we want to examine the scalability with respect to the number of subdomains once we substitute the direct solution method. Secondly, we want to observe what level of accuracy is needed at the subdomain level such that the outer number of iterations remains within a satisfactory range.

\FloatBarrier
\vspace*{-1.5mm}
\subparagraph{High-tolerance: $10^{-10}$}
We start with a tolerance of $10^{-10}$ with results given in Table \ref{bootland_mini_03_tab:02}. This case is the closest to the use of a direct solver (see Table \ref{bootland_mini_03_tab:01}). Comparing, we observe that the results are almost identical when using 4 subdomains. Once we increase the number of subdomains $N$, the number of FGMRES iterations increases for both 10 and 20 ppwl. However, the increase is more noticeable when using 10 ppwl. For example, when $k = 160$ and $N = 9$ a direct solver on the subdomains leads to 262 iterations while the inexact approach converges in 301 outer iterations. However, when we double $n_\text{ppwl}$ to 20, we go from 289 to 292 outer iterations. In all cases, as expected, the number of outer iterations increases with the wave number $k$. A similar yet slower increase in iteration counts is observed for the average number of inner iterations required by the deflated GMRES approach on the subdomains. As mentioned previously, the deflation preconditioner becomes more efficient on finer meshes. This can also be observed in our results: while the number of outer FGMRES iterations (using RAS) increases when moving from 10 ppwl to 20 ppwl, the number of inner GMRES iterations (using two-level deflation) decreases as the local subdomain systems become larger. Additionally, in this case, the number of inner iterations appears to be scaling better with the wave number $k$.

\begin{table}[t]
\caption{FGMRES iteration counts using the one-level RAS preconditioner with subdomain problems solved inexactly to a relative tolerance of $10^{-10}$ using GMRES with a two-level deflation preconditioner. In parentheses we display the average number of GMRES iterations per subdomain solve.}
\begin{minipage}{.49\linewidth}
\centering
10 ppwl\\
\scalebox{1}{
\begin{tabular}{@{}cc|cccc@{}}
\hline
& & \multicolumn{4}{c}{$N$} \\ \hline \hline
$k$ & $n_{\text{glob}}$ & 4 & 9 & 16 & 25 \\ \hline
20 & 30 & 20 (23) & 27 (25) & 56 (24) & 46 (21) \\
40 & 60 & 32 (35) & 62 (29) & 86 (26) & 101 (25) \\
80 & 120 & 65 (45) & 137 (32) & 192 (32) & 221 (29) \\
160 & 240 & 160 (63) & 301 (58) & 373 (36) & 518 (33) \\ \hline
\end{tabular}}\label{bootland_mini_03_tab:02a}
\end{minipage}%
\hfill
\begin{minipage}{.49\linewidth}
\centering
20 ppwl\\
\scalebox{1}{
\begin{tabular}{@{}cc|cccc@{}}
\hline
& & \multicolumn{4}{c}{$N$} \\ \hline \hline
$k$ & $n_{\text{glob}}$ & 4 & 9 & 16 & 25 \\ \hline
20 & 60 & 20 (30) & 43 (26)& 42 (25) & 59 (23)\\
40 & 120 & 37 (31) & 66 (32) & 93 (27) & 112 (27)\\
80 & 240 & 75 (36)& 132 (30) & 191 (30)& 268 (27)\\
160 & 480 & 131 (53) & 292 (47) & 407 (31) & 530 (28) \\ \hline
\end{tabular}}\label{bootland_mini_03_tab:02b}
\end{minipage}\label{bootland_mini_03_tab:02}
\end{table}
\FloatBarrier

\vspace*{-1.5mm}
\subparagraph{Medium-tolerance: $10^{-5}$}
In Table \ref{bootland_mini_03_tab:03} we report the results when lowering the inner tolerance to $10^{-5}$. We compare with the results reported in Table \ref{bootland_mini_03_tab:02}. A general observation is that as $k$ increases, so does the number of outer FGMRES iterations. Naturally, lowering the inner tolerance ensures we require less iterations to converge on the subdomains.

For the largest wave number reported and 9 subdomains, we needed 262 outer iterations when using a direct solver. For the inexact approach with tolerance $10^{-5}$ the number of outer iterations increases to 301, which is the same as when using a tolerance of $10^{-10}$. However, for the finer mesh with 20 ppwl the number of outer iterations goes up from 292 to 308. At the same time, the number of inner iterations reduces accordingly: from 58 to 40 for 10 ppwl and from 47 to 31 for 20 ppwl.

If we increase the number of subdomains from 9 to 25, the outer number of FGMRES iterations increases more rapidly. If we use 20 ppwl, the direct local solves lead to 520 outer iterations. This goes up to 555 when we use the iterative approach and a tolerance of $10^{-5}$. Note that the extra outer iterations compared to a tolerance of $10^{-10}$ is surmountable as relaxing the tolerance by 5 orders of magnitude leads to an increase of 25 iterations (from 530 to 555). Moreover, if we compare the number of inner iterations, a finer mesh works better with two-level deflation preconditioned GMRES on the subdomains, since we now need 15 iterations on average.

Similarly, on the finer mesh the number of inner iterations scales better with increasing wave number $k$ by adding more subdomains. Contrary to the results for $k < 160$, moving from 10 ppwl to 20 ppwl with 25 subdomains leads to less outer iterations. Thus, for larger wave numbers, using a finer mesh with more subdomains leads to a smaller number of outer and inner iterations. This effect is not observed with respect to the direct solves on the subdomains and/or the use of the tolerance $10^{-10}$ (see Table \ref{bootland_mini_03_tab:01} and Table \ref{bootland_mini_03_tab:02}): here as we go from 10 ppwl to 20 ppwl, the number of outer iterations always increases.

\begin{table}[t]
\caption{FGMRES iteration counts using the one-level RAS preconditioner with subdomain problems solved inexactly to a relative tolerance of $10^{-5}$ using GMRES with a two-level deflation preconditioner. In parentheses we display the average number of GMRES iterations per subdomain solve.}
\begin{minipage}{.49\linewidth}
\centering
10 ppwl\\
\scalebox{1}{
\begin{tabular}{@{}cc|cccc@{}}
\hline
& & \multicolumn{4}{c}{$N$} \\ \hline \hline
$k$ & $n_{\text{glob}}$ & 4 & 9 & 16 & 25 \\ \hline
20 & 30 & 20 (12) & 27 (14) & 65 (14) & 46 (12)\\
40 & 60 & 34 (20) & 76 (15)& 95 (14)& 122 (14) \\
80 & 120 & 72 (26) & 154 (17)& 210 (19)& 262 (17)\\
160 & 240 & 175 (44) & 301 (40) & 398 (19) & 572 (20) \\ \hline
\end{tabular}}\label{bootland_mini_03_tab:03a}
\end{minipage}%
\hfill
\begin{minipage}{.49\linewidth}
\centering
20 ppwl\\
\scalebox{1}{
\begin{tabular}{@{}cc|cccc@{}}
\hline
& & \multicolumn{4}{c}{$N$} \\ \hline \hline
$k$ & $n_{\text{glob}}$ & 4 & 9 & 16 & 25 \\ \hline
20 & 60 & 20 (15) & 50 (13) & 43 (13) & 59 (12) \\
40 & 120 & 37 (16) & 76 (18)& 104 (14)& 118 (14)\\
80 & 240 & 86 (19)& 148 (16)& 218 (16)& 314 (14)\\
160 & 480 & 144 (34) & 308 (31) & 431 (16) & 555 (15) \\ \hline
\end{tabular}}\label{bootland_mini_03_tab:03b}
\end{minipage}\label{bootland_mini_03_tab:03}
\end{table}

\subparagraph{Low-tolerance: $10^{-2}$}
Finally, we reduce the inner tolerance to just $10^{-2}$ and report results in Table \ref{bootland_mini_03_tab:04}. While the overall observations follow a similar trend to the previous case, the number of inner iterations are reduced drastically. This comes at the expense of a higher number of outer iterations as the wave number and number of subdomains increase.

The most noticeable result is again for the highest wave number, $k = 160$. If we use 20 ppwl, the direct local solves lead to 520 outer iterations. This goes up to 584 when we use the iterative approach with a tolerance of $10^{-2}$. The extra outer iterations compared to a tolerance of $10^{-10}$ is again surmountable as relaxing the tolerance by 8 orders of magnitude leads to an increase of 54 outer iterations (from 530 to 584). Meanwhile, the average inner iterations goes from 28 (for $10^{-10}$) to 15 (for $10^{-5}$), and finally to 6 iterations when using a tolerance of $10^{-2}$.

Analogous to the case where we set the tolerance to $10^{-5}$, we again observe that, as an exception to the rule that increasing $n_\text{ppwl}$ leads to more outer iterations, the number of outer iterations actually decreases when using 20 ppwl instead of 10 ppwl. This effect is only observed for the iterative approach on the subdomains in combination with a sufficiently low tolerance, here $10^{-5}$ or $10^{-2}$.

An important take-away message here is that the outer iteration of the one-level RAS preconditioned FGMRES method is able to reach convergence even when the subdomain systems are solved only to a relatively low level of accuracy.
\begin{table}[t]
\caption{FGMRES iteration counts using the one-level RAS preconditioner with subdomain problems solved inexactly to a relative tolerance of $10^{-2}$ using GMRES with a two-level deflation preconditioner. In parentheses we display the average number of GMRES iterations per subdomain solve.}
\begin{minipage}{.49\linewidth}
\centering
10 ppwl\\
\scalebox{1}{
\begin{tabular}{@{}cc|cccc@{}}
\hline
& & \multicolumn{4}{c}{$N$} \\ \hline \hline
$k$ & $n_{\text{glob}}$ & 4 & 9 & 16 & 25 \\ \hline
20 & 30  & 20  (5)  & 27  (6) & 73  (7) & 50   (5)\\
40 & 60  & 42  (8)  & 87  (7) & 109 (7) & 126  (6)\\
80 & 120 & 84  (13) & 172 (8) & 241 (10)& 298  (8)\\
160& 240 & 211 (30) & 332 (22)& 451 (9) & 1007 (8)\\ \hline
\end{tabular}}\label{bootland_mini_03_tab:04a}
\end{minipage}%
\hfill
\begin{minipage}{.49\linewidth}
\centering
20 ppwl\\
\scalebox{1}{
\begin{tabular}{@{}cc|cccc@{}}
\hline
& & \multicolumn{4}{c}{$N$} \\ \hline \hline
$k$ & $n_{\text{glob}}$ & 4 & 9 & 16 & 25 \\ \hline
20 & 60  & 21  (7) & 59  (5) & 49  (5) & 72  (5)\\
40 & 120 & 44  (7) & 84  (8) & 124 (6) & 134 (6)\\
80 & 240 & 94  (9) & 154 (7) & 229 (8) & 333 (6)\\
160& 480 & 154 (20)& 327 (19)& 450 (7) & 584 (6)\\ \hline
\end{tabular}}\label{bootland_mini_03_tab:04b}
\end{minipage}\label{bootland_mini_03_tab:04}
\end{table}

To provide some perspective on these results, we repeat the 10 ppwl experiment but now use GMRES preconditioned by ILU(0) as the subdomain solution method (with tolerance $10^{-2}$). The results in Table \ref{bootland_mini_03_tab:05} show that the average number of inner iterations drastically increases. Further, for $k = 160$ simulation run times were noticeably increased. Note that, for $k = 160$ and $N = 25$, using the two-level deflation preconditioner with 20 ppwl leads to both a lower inner and outer iteration count.
\begin{table}[t]
\caption{FGMRES iteration counts using the one-level RAS preconditioner with subdomain problems solved inexactly to a relative tolerance of $10^{-2}$ using GMRES with an ILU(0) preconditioner. In parentheses we display the average number of GMRES iterations per subdomain solve. Here we use 10 ppwl.}
\centering
\scalebox{1}{
\begin{tabular}{@{}cc|cccc@{}}
\hline
& & \multicolumn{4}{c}{$N$} \\ \hline \hline
$k$ & $n_{\text{glob}}$ & 4 & 9 & 16 & 25 \\ \hline
20  & 30  & 28  (21)  & 42  (12)  & 74  (10) & 53  (7)  \\
40  & 60  & 44  (64)  & 80  (36)  & 110 (25) & 131 (17) \\
80  & 120 & 89  (250) & 177 (121) & 239 (83) & 303 (55) \\
160 & 240 & 221 (983) & 344 (502) & 475 (260)& 658 (199)\\ \hline
\end{tabular}}\label{bootland_mini_03_tab:05}
\end{table}
\FloatBarrier

\section{Conclusions}
\label{bootland_mini_03_sec:4}

In this work we examined the utility of the one-level RAS preconditioner together with FGMRES to solve the 2D homogeneous Helmholtz equation when using an inexact solution method for the subdomain solves. Our results support the notion that the direct solve can be substituted by an efficient iterative solver. By using two-level deflation as a local preconditioner, we are able to keep the number of inner iterations on the subdomains low and scalable with respect to the wave number $k$.

The next step would be to include a coarse space and experiment with a two-level RAS preconditioner combined with inexact solves on the subdomains. Adding a coarse space would reduce the impact on the number of outer iterations when substituting direct solves for inexact solves on the subdomains. The trade-off between a higher number of outer iterations and a fast and memory efficient local subdomain solve needs to be analysed in large-scale applications to determine the break-even point in terms of wall-clock time, identifying where the iterative approach can be beneficial. Especially in high-frequency 3D applications, the inclusion of the coarse space can become a bottleneck. To reduce the outer iteration count, we either need to solve with a large coarse space or on larger subdomains. Both options can be costly when using a direct method and so an inexact solver is likely more suitable. Further, in the iterative approach the inner Krylov solvers for the subdomain problems may also benefit from the use of recycling techniques, which could further reduce the number of inner iterations and increase efficiency.

\bibliographystyle{spmpsci}
\bibliography{bootland_mini_03}

\end{document}